\documentclass{article}
\usepackage{amssymb}
\title{A Useful Property of the Finite Nonabelian Groups}
\author{Leendert Bleijenga}
\date{11 March 2014}
\begin{document}
\maketitle
\begin{abstract}
In version $v1$ (under a different title) I was trying to give a new proof of Wedderburn's Little Theorem $(WLT)$, stating that a finite division ring is commutative, but I failed. So I had to withdraw the paper (version $v2$). Firstly I became aware of a new theorem. So in the mean time in version $v3$ under a new title (see top of this page) I proved a useful property of the finite nonabelian groups stating that one of its maximally abelian subgroups is not an eigenheimer. In the present version $v4$ I prove finally $WLT$ in a new section $2$ entitled: A New Proof of Wedderburn's Little Theorem: A Finite Division Ring is Commutative. 
\end{abstract}

\setlength{\parskip}{5mm}

Subjects: Group Theory; Rings and Algebras;


\footnote{Tobbe of geen tobbe: dat was de was.}

\footnote{Email: leen.bleijenga@gmail.com}

\newpage

\section{Eigenheimers}

We begin with some facts about finite nonabelian groups.

\newtheorem{defi}{Definition}
\begin{defi}
Let $G$ be a finite nonabelian group with center $Z$ and let the $r$ maximally abelian subgroups be $H_{1}, H_{2}, \cdots, H_{r}, r \geq 3$. $G$ is called $Z$-indepen\-dent if $H_{i} \cap H_{j} = Z$ for all $i, j, i \not = j$ and $G$ is called $Z$-dependent if there exists $i, j$ $i \not = j$ and $H_{i} \cap H_{j} = D \supset Z$ and $D \not = Z$.
\end{defi}

\newtheorem{theo1}[defi]{Theorem}
\begin{theo1}
Let $G$ be a finite nonabelian group with center $Z$ and let $\mathcal{H} = \{H_{1}, H_{2}, \cdots, H_{r}\}$ be the set of the $r$ maximally abelian subgroups of $G$. Then the following $6$ claims hold: $(i)$ $Z \subset H_{i}, i = 1, 2, \cdots, r$. $(ii)$ $\cup_{i = 1}^{i = r} H_{i} = G$. $(iii)$ $r \geq 3$. $(iv)$ $\cap_{i = 1}^{i = r} H_{i} = Z, r \geq 3$. $(v)$ Let $G$ be $Z$-dependent and suppose that $\mathcal{H}_{D} = \{H'_{1}, H'_{2}, \cdots, H'_{s}\}$ is the subset of all maximally abelian subgroups of $G$ which contain subgroup $D \supset Z$, $D \not = Z$ and let $H'_{1} \cap H'_{2} = D$. Then the subgroup 
$<\mathcal{H}_{D}> = H'_{1} \cup H'_{2} \cup \cdots \cup H'_{s} < G$. Moreover, the maximally abelian subgroups of $<\mathcal{H}_{D}>$ are those of $G$ containing $D$ and are thus in $\mathcal{H}_{D}$. $(vi)$ If $\mathcal{H} = \mathcal{H}_{1} \cup \mathcal{H}_{2}$ is a partition of $\mathcal{H}$ and $\mathcal{H}_{1} \cap \mathcal{H}_{2} = \emptyset$ then $<\mathcal{H}_{1}> = G$ or $<\mathcal{H}_{2}> = G$.
\end{theo1}

Proof: $(i)$ We prove first that $Z \subset H_{1} = H$. Consider the set $HZ$. Let $h_{1}, h_{2} \in H$ and $z_{1}, z_{2} \in Z$. Then $h_{1} z_{1} \times h_{2} z_{2} = h_{1} h_{2} z_{1} z_{2} = h_{3} z_{3} \in HZ$ and $(h_{1} z_{1})^{-1} = z_{1}^{-1} h_{1}^{-1} = h_{1}^{-1} z_{1}^{-1} \in HZ$ and $HZ$ is a subgroup of $G$. But $HZ$ is also abelian for $h_{1} z_{1} \times h_{2} z_{2} = h_{2} z_{2} \times h_{1} z_{1}$. Now $H \subset HZ$ but $H$ is a maximally abelian subgroup of $G$. Thus $H = HZ$ and $Z \subset HZ = H$. Thus all maximally abelian subgroups contain $Z$.

$(ii)$ Every element $x \in G$ generates a cyclic subgroup (which is abelian) and lies in one or more maximally abelian subgroups.

$(iii)$ If there is only one maximally abelian subgroup $H_{1}$ then $H_{1} = G$. A contradiction. And if there are two maximally abelian subgroups $H_{1}$ and $H_{2}$ then $H_{1} \cup H_{2} = G$ and $G$ is the union of two proper subgroups which is impossible as we will demonstrate. Let $x \in H_{1} - H_{2}$ and $y \in H_{2}- H_{1}$. Then, for example, $xy \in H_{1}$ and $x \in H_{1}$ and it follows that $y \in H_{1}$. A contradiction.

$(iv)$ Let $D = \cap_{i = 1}^{i = r} H_{i}$ so that $D \supset Z$ is abelian. Let $D \not = Z$ and let $d \in D - Z$ Then $d$ commutes with all the elements in $H_{1}, H_{2}, \cdots, H_{r}$. Thus $d$ commutes with all the elements in $\cup_{i = 1}^{i = r} H_{i} = G$ and $d \in Z$. A contradiction. Thus $D = Z$.

$(v)$ $G \not= <\mathcal{H}_{D}>$ because otherwise $D \subset Z$. Let $x \in <\mathcal{H}_{D}>$ and $x \not \in H'_{i}, i = 1, 2, \cdots, s$. Now the elements of $D$, which is abelian, commute with all elements of $<\mathcal{H}_{D}>$ and thus also with $x$. $D$ and $x$ generate an abelian subgroup which lies in a maximally abelian subgroup of $G$, let us say $H'_{s + 1}$. Now $D \subset H'_{s + 1} \in\mathcal{H}_{D}$ contradicting the maximality of $s$. (The join is equal to the union.) Let $E$ be a maximally abelian subgroup of $<\mathcal{H}_{D}>$ thus $D \leq E$. Now $E$  lies in a maximally abelian subgroup $F$ of $G$ and $D \leq F$ so that $F \in \mathcal{H}_{D}$. Thus $E = F$. This completes the proof.

$(vi)$ Otherwise $G$ is the union of two proper subgroups, which is impossible.

We need a definition.

\begin{defi}
A subgroup $H$ of a group $G \neq H$ is called an \textbf{eigenheimer} of $G$ (or \textbf{self-normalizing}) if the normalizer of $H$ in $G$ is equal to $H$. Thus $N_{G}(H) = H$.
\end{defi}

\begin{theo1}
  Let $G$ be a finite group with eigenheimer $H$ and let $x \in G$. Then the conjugate subgroup $xHx^{-1}$ is also an eigenheimer of $G$.
\end{theo1}

Proof: If $xHx^{-1}$ is not an eigenheimer than $xHx^{-1} < N_{G} (xHx^{-1}) = x N_{G} (H) {x}^{-1}$ so that $H < N_{G}(H)$. A contradiction. Let's prove the last equality: $(\subset)$ Let $a \in N_{G} (xHx^{-1})$ then $axHx^{-1}a^{-1} = xHx^{-1}$ or $x^{-1}ax H x^{-1}a^{-1}x = H$ so that $x^{-1}ax \in N_{G} (H)$ or $a \in xN_{G} (H) x^{-1}$. We now prove $(\supset)$: Let $b \in x N_{G} (H) x^{-1}$ or $x^{-1}bx \in N_{G} (H)$ so that $x^{-1}bx Hx^{-1}b^{-1}x = H$ or $bxHx^{-1}b^{-1} = x H x^{-1}$ or $b \in N_{G} (xHx^{-1})$ and the equality follows.

We prove a few grouptheoretic theorems:

\begin{theo1}
  Let $G$ be a finite group and let $H$ be a proper subgroup of $G$ Then $\cup_{x \in G} (xHx^{-1}) \not = G$.
\end{theo1}

Proof: Let $|G| = g$, $|H| = h$. Let $D = \bigcap_{x \in G} (xHx^{-1})$ and let $|D| = d \geq 1$. Then $|G : N_{G} (H)| \leq \frac{g}{h}$ and thus $(g - d) \leq (\frac{g}{h}) (h - d)$ or $h(g - d) \leq g(h - d)$ or $gd \leq hd$ or $g = h$. A contradiction.

\newtheorem{theo4}[defi]{The Main Theorem for Finite Nonabelian Groups}
\begin{theo4}
  Let G be a finite nonabelian group with center $Z$. Then one of the maximally abelian subgroups of $G$ is not an eigenheimer of $G$.
\end{theo4}

Proof: We use induction on the number $|G|$ and we assume that $G$ is a minimal counterexample. There are two cases to be considered:

$(i)$ Firstly that $G$ is $Z$-independent.
Let $H$ and $K$ be maximally abelian subgroups of $G$ and let $H \not = K$. Then $H \cap K = Z$. We assume that all maximally abelian subgroups are eigenheimers, otherwise we are done. Let $|G| = g$, $|H| = h$, $|K| = k$ and $|Z| = z$. Suppose $(i a)$ that all maximally abelian subgroups are conjugate with $H$. Then $g - z = \frac{g}{h} (h - z)$ or $g = h$. A contradiction. Suppose $(i b)$ $H$ and $K$ are not conjugate. Then $(g - z) \geq \frac{g}{h} (h - z) + \frac{g}{k} (k - z)$ or without loss of generality $(g - z) \geq 2 \frac{g}{h} (h - z)$ or $h (g - z) \geq 2 g (h - z)$ or $0 \geq hz + g (h - 2z) \geq hz$. A contradiction. From these two contradictions it follows that at least one of the maximally abelian subgroups, let us say $H$, is not an eigenheimer: $N_{G} (H) \neq H$.

$(ii)$ Secondly that $G$ is Z-dependent. Referring to theorem $2$ $(v)$ and knowing that the subgroup $<\mathcal{H}_{D}> \not = G$
we use induction on $|<\mathcal{H_{D}}>|$. As the maximally abelian subgroups of $<H_{D}>$ are maximally abelian subgroups of $G$, then by induction one of the $H'_{i}, i = 1, 2, \cdots, s$ is not an eigenheimer in $<\mathcal{H_{D}}>$ and thus also not in $G$. This completes the proof.

\begin{theo1}
Let $G$ be a finite nonabelian group such that the proper subgroups of $G$ are all abelian. Then one of the maximal subgroups of $G$ is an abelian normal subgroup of $G$.
\end{theo1}

Proof: One of the maximal subgroups, which is abelian, is not an eigenheimer and thus a normal subgroup of $G$.

\newtheorem{theo3}[defi]{Theorem of Zassenhaus [1952]}
\begin{theo3}
Let $G$ be a finite group and for every abeli\-an subgroup $H$ of $G$ we have $N_{G} (H) = C_{G} (H)$. Then $G$ is abelian.
\end{theo3}

Proof: Suppose that $G$ is nonabelian. Let $H_{i}$ be the maximally abelian subgroups of $G$. Then one of the maximally abelian subgroups, let us say $H$, is not an eigenheimer by theorem $6$. Thus $H < N_{G} (H) = C_{G} (H) \neq H$. Let $c \in C_{G} (H) - H$ then the group generated by $H$ and $c$ is abelian. A contradiction. Thus $G$ is abelian.

\section{A New Proof of Wedderburn's Little Theorem: A Finite Division Ring is Commutative}

(This was also the original title of version $v1$, so I'm happy again.)

We will prove the well-known fact that a finite division ring is commutative (Wedderburn, 1905) along the following lines. We start with a minimal counterexample $L$, with center $Z$, viz. a finite division ring which is not commutative but all its maximal division subrings are commutative and are thus subfields. We say that $L$ is $Z$-dependent or $Z$-independent when we mean actually that $L^{\times}$ is $Z^{\times}$-dependent or $Z^{\times}$-independent. There are two possibilities to be considered. The first case $(I)$ is where $L$ is $Z$-dependent and the second case $(II)$ is where $L$ is $Z$-independent, see Definition $1$.

Case $(I)$ Let $\mathcal{F} = \{F_{1}, F_{2}, \cdots, F_{r}\}$ be the set of the $r$ maximal subfields of the finite division ring $L$. All the maximal subfields contain the center $Z$ of $L$ and let $F_{1} \cap F_{2} = D \supset Z$, and $D \not = Z$. Thus $L$ is $Z$-dependent. $D$ is also a commutative subfield of $L$ but is not a maximal subfield. Let $\mathcal{F}'$ = $\{F'_{1}, F'_{2}, \cdots, F'_{s}\}$ be the set of maximal subfields containing also $D$. Then $\bigcup \mathcal{F}'$ is not only a multiplicative group (if we drop the zero) (see Theorem 2 (v)) but also an additive subgroup (if we allow the zero) and thus a proper division ring which is by induction commutative as we now will show the additive case. Let $a \in F'_{i}$ and $b \in F'_{j}$ then $a^{-1}b$ is an element of let us say $F'_{k} \in \mathcal{F}'$. But $1 - a^{-1}b$ is also in $F'_{k}$. Thus the product $a(1 - a^{-1}b) = a - b$ is in let us say $F'_\ell \in \mathcal{F}'$. We have a contradiction.

Case $(II)$ Let $\mathcal{F} = \{F_{1}, F_{2}, \cdots, F_{r}\}$ be the set of the $r$ maximal subfields of the finite division ring $L$. All the maximal subfields contain the center $Z$ of $L$. We assume now that $F_{i} \cap F_{j} = Z$ for all $i$ and $j$, $i \not = j$ and $1 \leq i, j \leq r$. Thus $L$ is $Z$-independent. According to main theorem $6$ one of the maximally abelian subgroups $F^{\times}_{i}, i = 1, \cdots, r$ is not an eigenheimer. Let us say that $F^{\times}_{1}$ is not an eigenheimer so that $N:= N_{L^{\times}} (F^{\times}_{1}) \not = F^{\times}_{1}$. We assume that $|F_{1}| = p^{m}$ and $|Z| = p^{z}$ so that $z|m$. Let $zs|m$ be such that $s$ is a prime number then  $F_{1}$ contains a subfield $E$ such that $|E| = p^{zs}$. But then also $N = N_{L^{\times}} (E^{\times})$, as we shall now prove. Let us recall that a finite field minus the zero is a cyclic group, a subgroup of a finite cyclic group is also cyclic and that two subgroups of a finite cyclic group with the same order are equal and moreover that the quotient group of a finite cyclic group with one of its subgroups is also cyclic. Suppose that $N_{1}$ is the normalizer of $F_{1} - \{0\}$ and $N_{2}$ is the normalizer of $E - \{0\}$. Let $x \in N_{1}$ then $x F_{1} x^{-1} =F_{1}$ which implies by the aforementioned arguments that $x E x^{-1} = E$. Thus $x \in N_{2}$. Now let $y \in N_{2}$ and thus $y E y^{-1} = E$. Suppose that $y F_{1} y^{-1} \neq F_{1}$ then $y F_{1} y^{-1} \cap F_{1} = E$ but $y F_{1} y^{-1} \cap F_{1} = Z$ because the two fields are $Z$-independent. Thus also $y F_{1} y^{-1} = F_{1}$ and $y \in N_{1}$ and $N_{1} = N_{2} = N$ what we had to prove. Now let $x \in N - F_{1}$ and we assume that $x$ is also an element of $F_{2}$. Let $d := x^{k} \in Z^{\times}$, $k$ a minimal positive integer. Consider the left vector space $V = F_{1} + F_{1} x^{1} + \cdots + F_{1} x^{k - 1}$ and the left vector space $W = E + E x^{1} + \cdots + E x^{k - 1}$. We will prove that $V = W = L$ by showing that $V$ and $W$ are division rings. Let $v_{1} = (\cdots + f_{i} x^{i} + \cdots) \in V$
and $v_{2} = (\cdots + f_{j} x^{j} + \cdots) \in V$ then $v_{1}v_{2} = (\cdots + f_{i} x^{i} f_{j} x^{j} + \cdots) = (\cdots + f_{i} (x^{i} f_{j} x^{-i}) x^{i + j} +\cdots) \in V$, for if $i + j \geq k$ then $x^{i + j} = dx^{i + j - k}$. Thus $V$ is a division ring and $V = L$. In the same way we can prove that $W$ is also a division ring and thus $W = L = V$. We saw earlier  that $B = \{1, x, x^{2}. \cdots, x^{k - 1}\}$ is a basis of $Z(x) \subset F_{2}$ and we prove now that it is a basis of $L$ over $F_{1}$. Suppose that a proper subset $B_{1} = \{x_{1}, x_{2}, \cdots, x_{m}\}$ of $B$ is a basis of $L$ over $F_{1}$. Then for $x_{m + 1} \in B - B_{1}$ we can write: $x_{m + 1} = f_{1} x_{1} + f_{2} x_{2} + \cdots + f_{m} x_{m}$. Then $x x_{m + 1} x^{-1} = x_{m + 1} = xf_{1}x^{-1} x_{1} + xf_{2}x^{-1} x_{2} + \cdots + xf_{m}x^{-1} x^{m}$. So $f_{j} = xf_{j}x^{-1}$ and $f_{j} \in Z, j = 1, 2, \cdots, m$. From this contradiction we have $E = F_{1}$ and $|L| = p^{zs \times k}$. 
Let us define the group $C := <x>$ and group $D := <x^{k}>$. Then $C$ is a subgroup of $F_{2}^{\times}$ and $D$ is a subgroup of $Z^{\times}$. Let $G := C/D = \{D, xD, x^{2}D, \cdots, x^{k - 1}D\}$ which is a cyclic group of order $k$. We construct a cyclic group $H = \{\phi, \phi^{2}, \cdots, \phi^{k} \}$ of order $k$ of automorphisms $\phi^{i}$ of the field $F_{1}$ in the following way. With the element $\phi^{i} \in H$ we define $\phi^{i} (a):= x^{i} a x^{-i}$, $a \in F_{1}$, then $\phi^{i}$ is an automorphism of the field $F_{1}$ which leaves the elements of $Z$ invariant. Moreover $\phi^{k} = 1$ for $\phi^{k} (a) = x^{k} a x^{-k} = dad^{-1} = a, a \in F_{1}$. We see also that $\phi^{i + 1} = \phi(\phi^{i})$ and that $\phi^{i} (a + b) = x^{i} (a + b) x^{-i} = x^{i} a x^{-i} + x^{i} b x^{-i} = \phi^{i} (a) + \phi^{i} (b)$ and that $\phi^{i} (ab) = x^{i} (ab) x^{-1} = (x^{i} a x^{-i}) (x^{i} b x^{-i}) = \phi^{i} (a) \phi^{i} (b),  a, b \in F_{1}$.
But according to a famous theorem of Artin which states that if $n$ automorphisms of a field $F$ form a group and $I$ is the subfield of the invariant elements under all of the $n$ automorphisms then $[F : I] = n$. So $|G| = |H| = [F_{1} : Z]$ or $k = s$. Then $|L| = p^{zs \times s}$. This means that all maximal subfields have the same number of elements, namely $|F_{i}| = p^{zs}$. The class equation sounds now: $L^{\times} = Z^{\times} + \sum (\frac{L^{\times}}{C(x)})$ where the summation is taken place over representatives of the remaining conjugacy classes which have by the way the same number of elements. So $p^{zs \times s} - 1 = p^{z} - 1 + \alpha (\frac{p^{zs \times s} - 1}{p^{zs} - 1})$. Multiplying the last equation by $p^{zs} - 1$ we see that $(p^{zs} - 1)(p^{z} - 1)$ is divisible by $p^{zs \times s} - 1$ but $p^{zs \times s} -1$ is divisible by $p^{zs} - 1$ leaving a quotient that is of the form $1 + p^{zs} + \cdots$ which is greater than $p^{z} - 1$. A contradiction. From the two contradictions in case $(I)$ and case $(II)$ it follows that every finite division ring is commutative. 

\textbf{References}
\begin{enumerate}
\item[1] Dr. Emil Artin, \textit{Galois Theory}, Edited and supplemented with a Section on Applications by Dr. Arthur N. Milgram, Second Edition, 1944, University of Notre Dame, Notre Dame, London. I found it on the Internet.
\item[2] A. Adrian Albert, \textit{Fundamental Concepts of Higher Algebra}, 1956, The University of Chicago Press.
\item[3] Nicolas Bourbaki, Elements of Mathematics, \textit{Algebra I}, Chapters 1-3. Softcover edition of the 2nd printing, 1989, \-Sprin\-ger.
\item[4] Nicolas Bourbaki, Elements of Mathematics, \textit{Algebra II}, Chapters 4-7, Translated by P.M. Cohn and J. Howie, Softcover printing of the 1st English edition of 1990, 2003, \-Sprin\-ger.
\item[5] M.I. Kargapolov and Ju.I. Merzljakov, \textit{Fundamentals of the Theory of \-Groups}, Translation from the Second Russian Edition by Robert G. Burns, 1979, \-Sprin\-ger.
\item[6] Prof.dr. J.H. van Lint and Dr. J.W. Nienhuys. \textit{Discrete Wiskunde}, 1991, Academic Service. (in Dutch)
\item[7] Dr. F. Loonstra, \textit{Inleiding tot de algebra}, zesde druk, 1979, Wolters-Noord\-hoff, Groningen. (in Dutch)
\item[8] Michio Suzuki, \textit{Group Theory I}, 1982, \-Sprin\-ger.
\item[9] Keith Conrad, \textit{Subgroups of Cyclic Groups}, http://www.math.uconn.edu/ ... I found it on the Internet.
\item[10] Anand Deopurkar, \textit{Subgroups of a finite cyclic group}, math.columbia.edu/ ... I found it on the Internet.
\item[11] ProofWiki, \textit{Quotient Group of Cyclic Group}, www.proofwiki.org/ ... I found it on the Internet.
\item[12]  Louis Weisner, Groups in which the normaliser of every element except identity is abelian, Bull. Amer. Math. Soc. 31 (1925), 413-416. I found it on Internet.
\item[13] Emil Artin, \"{U}ber einen Satz von Herrn J.H. Macgalan Wedderburn, Hamb. Abb. 5 (1928), pp. 245-250. I found it in EMIL ARTIN COLLECTED PAPERS, Edited by S. Lang and J.T. Tate, Springer-Verlag, 1965, pp. 301-306.
\item[14] H.E. Goheen, On a theorem of Zassenhaus, ams. org. 1954, pp. 799-800. I found it on the Internet.
\item[15] Shamil Asgarli, Wedderburn's Little Theorem. I found it on the Internet.  
\end{enumerate}
\end{document}